\documentstyle[12pt]{article}

\font\tenmsb=msbm10
\font\sevenmsb=msbm7
\font\fivemsb=msbm5

\newfam\msbfam
\textfont\msbfam=\tenmsb
\scriptfont\msbfam=\sevenmsb
\scriptscriptfont\msbfam=\fivemsb

\def\Bbb#1{{\fam\msbfam #1}}

\makeatletter
\ifnum\@ptsize=0  \fi \ifnum\@ptsize=2  \fi 

\newcommand\qed{{\hspace*{\fill}Q.E.D.\vskip12pt plus 1pt}}

\newcommand\sE{{\cal E}}
\newcommand\sG{{\cal G}}

\newcommand\sL{{\cal L}}
\newcommand\sO{{\cal O}}
\newcommand\sD{{\cal D}}
\newcommand\sC{{\cal C}}

\newcommand\zed{{\Bbb Z}}
\newcommand\comp{{\Bbb C}}
\newcommand\bP{{\Bbb P}}

\newcommand\proof{{\noindent\bf Proof.\ }}
\newtheorem{theorem}{Theorem}[section]
\newtheorem{lemma}[theorem]{Lemma}
\newtheorem{corollary}[theorem]{Corollary}
\newtheorem{proposition}[theorem]{Proposition}

\newtheorem{re}[theorem]{Remark}

\newtheorem{problem}[theorem]{Problem}

\begin{document}
\title {\bf Generalized Tsen theorem and \\
rationally connected Fano fibrations} \author{\bf Frederic
Campana, Thomas Peternell and Aleksandr Pukhlikov} \date{}

\maketitle



\parshape=1
3cm 10cm \noindent {\small \quad\quad\quad \quad\quad\quad\quad
\quad\quad\quad {\bf }\newline We prove that a fibration $X/{\bP}^1$ the general fiber of which is a smooth Fano three-fold is
rationally connected. The proof is based upon a generalization of
the classical Tsen theorem: we prove that a fibration $X/C$ over
a curve $C$ the general fiber of which is a Fano complete
intersection in a product of (weighted) projective spaces has a
section.} \vspace{1cm}

\vspace*{-0.5in}

\section{Introduction}


This paper is concerned with the following problem

\begin{problem}
Let $\pi\colon X\to {\bP}^1$ be a fibration, the general fiber
of which is rationally connected. Then $X$ is rationally
connected.

\end{problem}

A positive solution of this problem has the important consequence
that fibrations over rationally connected varieties with
rationally connected fiber are again rationally connected. If one
accepts the predictions of Mori theory, then one is reduced to
the analysis of Fano fibrations over ${\bP}^1$, however one
should have in mind that the fibers may all be singular in this
approach since the total space is singular in general when
passing from the original $X$ to a model admitting a Fano
fibration.

We have two main results:

\subsection{Generalized Tsen theorem}
Let ${\bP}$ denote the  weighted projective space ${\bP}(a_0,a_1,\dots,a_m)$. Let ${\cal V}$ denote the class of Fano
complete intersections in the product of weighted projective
spaces ${\bP}^{(1)}\times\dots\times {\bP}^{(k)}={\bP}^{\times}$
where ${\bP}^{(i)}={\bP}(a_0^{(i)},a_1^{(i)},\dots,a_{m(i)}^{(i)})$

\begin{theorem}{\bf  (the generalized Tsen theorem)} Let
$\pi\colon X\to C$ be a fibration over a curve, the general fiber
of which belongs to the class ${\cal V}$. Then there
exists a section $\beta\colon C\to X$, $\pi\circ\beta=
\mathop{\rm id}\nolimits_C$.
\end{theorem}

\begin{corollary} Let
$\pi\colon X\to {\bP}^1$ be a fibration over a rational curve,
the general fiber of which is a Fano complete intersection in the
product of weighted projective spaces. Then $X$ is rationally
connected.
\end{corollary}

\subsection{Fano fibrations}
\begin{theorem}  Let $\pi\colon X\to {\bP}^1$ be a
fibration, the general fiber of which is a smooth Fano
threefold. Then $X$ is rationally connected.
\end{theorem}

The proof of this theorem is by running through the classification
of Fano $3-$folds by Iskovskikh,Mori and Mukai and using ad hoc
methods in the single cases. We are giving details only in case $b_2 = 1$
while details for $b_2 \geq 2$ will be written up at a
different place.

After having completed our work up to some redigations, we got
aware of the new prepint [GHS01] in which the problem 1.1 appears to be completely solved with
totally different methods
from deformation theory and moduli theory of curves in the meantime also a preprint of de Jong-Starr appeared.
However this
paper remains meaningful, because the methods presented here are
effective (sections are explicitly constructed) and they
seem to be much simpler. An interesting aspect is also that our approach avoids the
involved considerations of multiple components of fibers.

\subsection{Plan of the proof of the generalized Tsen theorem}

We will prove a particular case of Theorem 2, when ${\bP}^{(i)}
\equiv{\bP}={\bP}^M$ is the ordinary projective space and the
general fiber $V\subset {\bP}\times\dots\times{\bP}$ is a
hypersurface of the type $(M,\dots,M)$. In the general case no
modification of the method is necessary, all the technique
developed below works automatically as it is (up to the Lefschetz
theorem 1.5). However, in the general case one needs considerably
more complicated notations. This is essentially the only reason
why we restrict ourselves by this typical particular case.

Our proof consists of two main steps:

{\begin{itemize}

\item
first, we construct a special model $X^{\sharp}/C$, birational to
$X/C$ over $C$:
$$
\begin{array}{ccc}
X & \stackrel{\chi}{-\,-\,\to}  & X^{\sharp} \\
\downarrow &  &  \downarrow \\
C & = & C,
\end{array}
$$
where over a point $x\in C$ of general position the birational map
$\chi$ is a {\it biregular} isomorphism of fibers, so that, in
particular, $\chi$ induces a bijection
$$
\mathop{\rm Sect}(X/C)\leftrightarrow\mathop{\rm
Sect}(X^{\sharp}/C)
$$
of the sets of sections and therefore the problems of existence of
a section of $X/C$ and $X^{\sharp}/C$ are equivalent.

Let us explain the meaning of this step. By assumption, a fiber
of general position of the fibration $X/C$ permits an embedding
$$
X_x=\pi^{-1}(x)\subset {\bP}\times\dots\times{\bP}
$$
as hypersurface of type $(M,\dots,M)$ in $\bP \times \ldots \times \bP = \bP^{\times k}$. However, in the
general case this embedding cannot be globalized, even on the
birational level, that is, there exists no birational map over $C$
$$
X-\,-\,\to{\bP}^{\times k}\times C
$$
to the trivial (or a locally trivial) ${\bP}^{\times
k}$-fibration over $C$, which would have embedded a fiber of
general position into ${\bP}^{\times k}$ as a hypersurface of
the type $(M,\dots,M)$. It is exactly for this reason that the
classical Tsen theorem [Ts36,Ko96] does not work. We will show, however
(Sec. 2), that there exists a certain analogue of such a
globalization: an embedding
$$
\begin{array}{ccc}
X^{\sharp} & \hookrightarrow  & {\bP}^{\sharp} \\
 &  &  \downarrow \\
 &  & C
\end{array}
$$
over $C$, where ${\bP}^{\sharp}$ is fibered over $C$ with the
general fiber ${\bP}^{\times k}$ (but of course ${\bP}^{\sharp}/C$ is not a locally trivial ${\bP}^{\times
k}$-fibration). The construction of the space ${\bP}^{\sharp}$
and the model $X^{\sharp}$ is very explicit and makes it possible
to study the problem of existence of sections.

\item
Second, we prove existence  of sections of
$X^{\sharp}/C\hookrightarrow{\bP}^{\sharp}/C$ (and
simultaneously give a description of all sections). We proceed in
two steps:

\begin{itemize}

\item
to begin with, we describe the sections of the ambient fibration
${\bP}^{\sharp}/C$.

\item
After that, using the explicit description of those divisors
$X^{\sharp}\subset{\bf P}^{\sharp}$, the fiber of which over a
point of general position $x\in C$ is a
$(M,\dots,M)$-hypersurface in ${\bf P}^{\times k}$, we make the
crucial step: we write out the conditions (a system of algebraic
equations on the parameters, defining the sections of the
fibration ${\bf P}^{\sharp}/C$), which make a section $s\colon
C\to{\bf P}^{\sharp}$ lie where it should:
$$
s(C)\subset X^{\sharp}.
$$
Now a simple dimension count in the spirit of the classical Tsen
theorem [Ts36],[Ko96] shows that these conditions really can be satisfied,
which completes the proof.

\end{itemize}

\end{itemize}}

There is one point to be a little careful concerning complete
intersections in weighted projective spaces:

\begin{proposition} Let
${\bP}^j$ be weigthed projective spaces, $ {\bP} = {\bP}^1
\times \ldots \times {\bP}^r $ and $X \subset {\bP}$ be a
quasi-smooth complete intersection of dimension at least $3$. Then the
restriction map
$$ {\rm Pic}({\bP}) \to {\rm Pic}(X)$$
is an isomorphism.
\end{proposition}

\proof Of course this is the classical Lefschetz theorem in case
${\bP}$ and $X$ are smooth.
In general there are three ways to conclude. \\
(a) Adapt the method of [Do81,3.2.4/3.2.5]
(where the case of one weighted projective space is treated). \\
(b) Use [Do80] (see remark (3.2.6) in [Do81]). \\
(c) Argue directly as follows. First we consider general quasi-smooth
complete intersections $X.$ Let ${\bP}_{m_j} \to {\bP}^j$ the
quotient map, let $\tilde {\bP} = {\bP}_{m_1} \times \ldots
\times {\bP}_{m_r} $ and consider the preimage $\tilde X$ of
$X$ in $\tilde {\bP}.$ Since $X$ is general, $\tilde X$ will be
smooth. Now apply the Lefschetz theorem for $\tilde X \subset
\tilde {\bP}$ and go back to $X.$ Of course this works only for
general $X.$ But any quasi-smooth $X$ is a deformation of those $X$ for
which we already know our claim. Then we see immediately that the
claim holds also for the special $X:$ it is of course clear that
${\rm Pic}(X) = {\rm Pic}(X_t)$, and a small argument shows that
every line bundle on $X$ is induced by a line bundle on ${\bP}$,
since it is topologically the case and the analytic cohomology in
the exponential sequence vanishes. \qed

\subsection{Acknowledgements}
For the third author, the present research has been carried out
in the framework of his Humboldt research project, and he is very
thankful to the Humboldt Foundation for the financial support.
The first author was supported by the Schwerpunktprogramm ``Global methods in complex geometry''
of the Deutsche Forschungsgemeinschaft.

\section{Construction of the special model}
In this section we realize the first step of our program: we
construct a special model of the fibration $X/C$, adjusted to our
program.

\subsection{The finite cover, parametrizing projections of the fibers}
Consider the fibration $X/C$ as in subsection 1.3. For a point $x\in C$ of general
position the fiber $\pi^{-1}(x)=X_x$ is a hypersurface of the type
$(M,\dots,M)$ in ${\bP}^{\times k}$. Therefore, for $X_x$ one
gets $k$ canonically given projections
$$
\mathop{\rm pr}\colon X_x\to{\bP}.
$$
These $k$ projections can be reconstructed from the variety $X_x$
itself and do not depend on its particular projective embedding.

\begin{proposition}  There exists a projective curve $B$
(reduced, but possibly reducible), which is a $k$-sheeted cover
of the curve $C$:
$$
\tau\colon B\to C,
$$
$\mathop{\rm deg}\tau=k$, the points of which parametrize the
projections of the fibers $X_x\to{\bP}$ in the following sense.
There exists a regular algebraic map
\begin{equation}
\label{1}
p\colon\pi^{-1}(C^{\circ})\mathop{\times}\nolimits_C
B\to{\bP},
\end{equation}
where $C^{\circ}\subset C$ is an open set of all such points $x\in
C$ that the fiber $X_x$ is a hypersurface of type
$(M,\dots,M)$ in ${\bP}^{\times k}$ (recall that by assumption
this is true for a fiber of general position only), the fiber
$\tau^{-1}(x)\subset B$ consists of precisely $k$ distinct points,
$$
\tau^{-1}(x)=\{y_1,\dots,y_k\},
$$
and the $k$ morphisms
$$
p|_{X_x\times\{y_i\}}\colon X_x\times\{y_i\}\to{\bP}
$$
are precisely $k$ projections
$$
\mathop{\rm pr}\colon X_x\to{\bP}\times\dots\times{\bP}\to{\bP}
$$
onto the factors for each point $x\in C^{\circ}$.

\end{proposition}

\proof  For simplicity let $k = 2.$ The fibers of the projections $X_x \to \bP^i,$ $x \in C, i = 1,2$ define
either two or one irreducible families of cycles. In the first case we obtain two rational maps $X \rightharpoonup \bP^i$,
which are regular on the general $X_x.$ Then we let $B$ be two disjoint copies of $C$ and the rest is clear. \\
So we may suppose that we have one irreducible family; let $\sC \to T$ denote its graph with projection $\rho: \sC \to X.$
Then we obtain canonically a map $h:T \to C$ which is of course a 2-sheeted cover (generically). Therefore we set $B = T.$
Over $C^0,$ we have $\sC = X \times_{C^0} B^0,$ and there exists a map $\sC \rightharpoonup \bP,$ just by sending
$z \in C_t$ to the point in $P$, to which $\rho(C_t)$ is mapped by the corresponding projection $X_h(z) \to \bP^i.$
\qed

\subsection{Construction of the ambient space ${\bP}^{\sharp}={\bP}^{B/C}$: an analytic
description over $C^{\circ}$}

Denote by the symbol ${\bP}_B={\bP}\times B$ the trivial
${\bP}$-fibration over $B$, and by $\pi_B\colon {\bP}_B\to B$ the corresponding projection.
Let ${\bP}_y={\bP}\times \{y\}= \pi^{-1}_B(y)$ be the fiber over a point $y\in
B$, and $\pi_C\colon{\bP}_B\to C$ be the composition
$\pi_C=\tau\circ \pi_B$. As above, denote by $C^{\circ}\subset C$,
the set of points of general position, which have precisely $k$
preimages in $B$.

An open (in the real topology) set $D\subset C^{\circ}$ (for
instance, a small disk) is said to be {\it simple}, if
$$
\tau^{-1}(D)=\coprod^k_{i=1}D_i,
$$
where $\tau\colon D_i\to D$ is an isomorphism for each $i$. Now
set for a simple $D\subset C^{\circ}$:
$$
{\bP}^{B/C}(D)=\prod^k_{i=1}{}_D({\bP}\times
D_i)\cong {\bP}^{\times k}\times D,
$$
where the symbol ${\prod}_D$ stands for the fiber product with
respect to the morphisms
$$
{\bP}\times D_i\to D_i\to D.
$$
Let
$$
\kappa_i\colon {\bP}^{B/C}(D)\to{\bP}
$$
mean the projection onto the $i$-th factor. Now set
$X(D)=\pi^{-1}(D)\subset X$ and look at the fiber product
$$
X(D)\mathop{\times}\nolimits_C B=X(D)\mathop{\times}\nolimits_C
\tau^{-1}(D) \simeq \coprod\limits^k_{i=1}X(D)^{(i)},
$$
where $X(D)^{(i)}$ is the identical copy of the set $X(D)$,
corresponding to the connected component $D_i$. Let
$$
\nu_i\colon X(D)\to X(D)^{(i)}
$$
be the corresponding isomorphism.

\begin{lemma} There exists a uniquely determined natural
morphism --- an inclusion
$$
\iota_D\colon X(D)\to {\bP}^{B/C}(D)
$$
such that for each $i=1,\dots,k$ the following diagram is
commutative:
\begin{equation}
\label{2}
\begin{array}{rcccl}
    & X(D) & \stackrel{\iota_D}{\longrightarrow} & {\bP}^{B/C}(D) &  \\
\nu_i & \downarrow &  & \downarrow & \kappa_i  \\
    &  X(D)^{(i)} & \stackrel{p_D}{\longrightarrow} & {\bP}, &
\end{array}
\end{equation}
where
\begin{equation}
\label{3} p_D\colon X(D)\mathop{\times}\nolimits_C B\to{\bP}
\end{equation}
is given by $p  \vert D \circ pr_D,$ with $p\vert D$ the restriction of the projection $p$ from Proposition 2.1
onto the set $D$.
\end{lemma}

\proof Define the morphism $\iota_D$ by the formula
$$
\iota_D=\prod^k_{i=1}(p\circ \nu_i)\times \pi\colon X(D)\to
{\bP}^{B/C}(D)=\left( \prod^k_{i=1}{\bP}\right)\times D.
$$
The diagram (\ref{2}) is now obviously commutative. It is also
obvious that (\ref{3}) determines the morphism $\iota_D$
uniquely. For a fiber $X_x$, $x\in D$, it means simply that the
$k$ projections
$$
p_y\colon X_x\to{\bP}_y,
$$
$y\in\tau^{-1}(x)$, determine the embedding
$$
\iota_x\colon X_x\to \prod_{y\in\tau^{-1}(x)}{\bP}_y.
$$
\qed

\begin{re} {\rm  Let $D^{\sharp}\subset D$ be an open subset. Then $D^{\sharp}$ is simple, too:
$$
\tau^{-1}(D^{\sharp})=\coprod^k_{i=1}D^{\sharp}_i,
$$
so that the construction, described above, can be applied to
$D^{\sharp}$, too. It is easy to see that there is a uniquely
determined natural inclusion
$$
\iota_{D^{\sharp},D}\colon{\bP}^{B/C}(D^{\sharp})\hookrightarrow{\bP}^{B/C}(D),
$$
which makes the following diagram commutative:
$$
\begin{array}{rcccl}
    & X(D^{\sharp}) & \hookrightarrow & X(D) &  \\
\iota_{D^{\sharp}} & \downarrow &  & \downarrow & \iota_D  \\
    &  {\bP}^{B/C}(D^{\sharp}) & \hookrightarrow & {\bP}^{B/C}(D). &
\end{array}
$$
This makes it possible to glue the varieties ${\bP}^{B/C}(D)$
and the inclusions $\iota_D\colon X(D)\hookrightarrow {\bP}^{B/C}(D)$
together and obtain an (open) variety ${\bP}^{B/C}(C^{\circ})$ and an embedding
$X(C^{\circ})\subset {\bP}^{B/C}(C^{\circ})$.}
\end{re}

\subsection{Construction of the ambient space ${\bP}^{B/C}$: an analytic
description in the ramification case}

To be able to work conveniently with the space ${\bP}^{B/C}$ it is
necessary to complete it over the points $x\in
C\setminus C^{\circ}$. We may assume that if the curve $B$ has
several irreducible components, then they are disjoint, that is, if a
fiber $\tau^{-1}(x)$ consists of less than $k$ points, then it
happens due to ramification of the morphism $\tau$. We also may assume
that the exceptional set $C\setminus C^{\circ}$ consists of such
points only, because if $\sharp \tau^{-1}(x)=k$, then the space
${\bP}^{B/C}$ is constructed exactly as above (although the
embedding morphism $i$ for the corresponding fiber
$X_x\hookrightarrow {\bP}^{B/C}_x$ is possibly no longer defined).

Now take a point $x\in C\setminus C^{\circ}$ and consider
$$
\tau^*x=\sum^a_{i=1}l_iy_i,
$$
where  $a\leq k-1$, $l_i\geq 1$ and at least one of the integers
$l_i$ is strictly bigger than 1. Let $D\subset C$ be a small disk
with the centre at the point $x$, where
$$
\tau^{-1}(D)=\coprod^a_{i=1}D_i, \quad D_i\ni y_i,
$$
and there exist local coordinates $z$ on $D$ and $z_i$ on $D_i$,
$z(x)=z_i(y_i)=0$, respectively, such that the morphism
$$
\tau\colon D_i\to D
$$
takes the form
$$
\tau\colon z_i\mapsto z_i^{l_i}.
$$
Let $\zeta_i$ be a primitive root of degree $l_i$ from 1 and
$G_i={\zed}_{l_i}$ be the cyclic group of order $l_i$.

Set
$$
Z(D)= \prod^a_{i=1}{}_D\left({\bP}^{\times l_i}\times
D_i\right)=
$$
$$
=\{((p^{(1)}_1;\dots;p^{(1)}_{l_1};z_1),
((p^{(2)}_1;\dots;p^{(2)}_{l_2};z_2),\dots,
((p^{(a)}_1;\dots;p^{(a)}_{l_a};z_a))|z^{l_1}_1=\dots=z^{l_a}_a\}.
$$

Now we define ${\bP}^{B/C}$ locally on $D$ in the following way:
$$
{\bP}^{B/C}(D)=Z(D)/G_1\times\cdots\times G_a,
$$
where the cyclic group $G_i=\langle\sigma_i\rangle$ acts on the
$i$-th factor of the fiber product $Z(D)$ as follows
$$
\sigma_i(p^{(i)}_1;\dots;p^{(i)}_{l_i};z_i)=
(p^{(i)}_{l_i};p^{(i)}_1;\dots;p^{(i)}_{l_i-1};\zeta_i z_i)
$$
(the points $p^{(i)}_j$ are permuted in the cyclic way). This
action is obviously compatible with the construction of the fiber
product ($(\zeta_i z_i)^{l_i}=z_i^{l_i}$) and thus the defined
action of $G_i$ can be naturally extended to an action on $Z^{\rm
an}(D)$ (on the other fiber factors the group $G_i$ acts
trivially).

It is checked easily that in the punctured disk
$D^{\sharp}=D\setminus\{x\}$ this construction coincides with
that of Sec. 2.2: more precisely, if $D^{\sharp}\subset D$ is a
simple open set, so that
$$
\tau^{-1}(D^{\sharp})=\coprod^k_{e=1}D^{\sharp}_e
$$
(in particular, $x\not\in D^{\sharp}$), then there is a uniquely
determined natural inclusion
$$
{\bP}^{B/C}(D^{\sharp})\hookrightarrow {\bP}^{B/C}(D),
$$
where on the left-hand side we have the space, constructed in
Sec. 2.2. Indeed,
$$
\tau^{-1}(D^{\sharp})\cap D_i=\coprod^{l_i}_{j=1}D^{\sharp}_{i,j},
$$
so that over $D^{\sharp}$ all the maps
$$
\tau_{i,j}\colon D^{\sharp}_{i,j}\to D^{\sharp}
$$
are isomorphisms. In particular, one can take $l_i$ single-valued
branches of the function $z^{1/l_i}=\tau^{-1}_{i,j}$. Fix one of
them and set
$$
\mu\colon D_i\cap\tau^{-1}(D^{\sharp})\to D^{\sharp}_i,
$$
$$
\mu\colon z_i\mapsto z^{\sharp}_i=(z_i^{l_i})^{1/l_i},
$$
where $D^{\sharp}_i=D^{\sharp}_{i,c(i)}$ is the subset
(connected component) $D^{\sharp}_{i,j}$, which corresponds to
the chosen branch. Now let
$$
\lambda\colon \prod^a_{i=1}{}_D \left({\bP}^{\times l_i}\times
(D_i\cap\tau^{-1}(D^{\sharp}))\right) = Z(D) \cap \left( \prod^a_{i=1} \bP^{\times l_i} \cap \tau^{-1} (D^\sharp)\right)
\to
\left(\prod^a_{i=1}\prod^{l_i}_{j=1}{\bP}\right)\times D
$$
be given on the $i$-th factor by the formula
$$
(p^{(i)}_1;\dots;p^{(i)}_{l_i};z_i)\mapsto
\sigma^f_i(p^{(i)}_1;\dots;p^{(i)}_{l_i};z_i),
$$
where $f\in{\zed}$ is defined by the condition
$$
z^{\sharp}_i=\zeta^f_i z_i.
$$
It is easy to check that $\lambda$ determines an isomorphism
$$
{\bP}^{B/C}(D)\cap \pi^{-1}(D^{\sharp})\cong {\bP}^{B/C}(D^{\sharp}),
$$
which, moreover, does not depend on the choice of $\zeta_i$ and
the single-valued branch of the root $(\cdot)^{1/l_i}$. This
isomorphism identifies ${\bP}^{B/C}(D^{\sharp})$ in the sense
of Sec. 2.2 with an open subset of the set ${\bP}^{B/C}(D)$.

What has been just said means that the open pieces ${\bP}^{B/C}(D)$, $D\subset C$,
glue together into a compact complex
space ${\bP}^{B/C}={\bP}^{B/C}(C)$ with a natural projection
$$
\pi\colon {\bP}^{B/C}\to C.
$$
Here for any point $x\in C^{\circ}$ the fiber of the projection
$\pi$ is
$$
{\bP}^{B/C}_x=\pi^{-1}(x)=\prod_{y\in\tau^{-1}(x)}{\bP}_y.
$$

\subsection{Construction of the ambient space ${\bP}^{B/C}$: an
algebraic description}

The analytic description of ${\bP}^{B/C}={\bP}^{B/C}_{\rm an}$
given above is very convenient for local computations. However, a
priori it is unclear whether ${\bP}^{B/C}$ is an algebraic
variety. So now we give a description of the space ${\bP}^{B/C}$
from the viewpoint of algebraic geometry, which we will
denote by ${\bP}^{B/C}={\bP}^{B/C}_{\rm alg}$; a little later
we will prove that there is an isomorphism of complex spaces
$$
{\bP}^{B/C}_{\rm an}\cong {\bP}^{B/C}_{\rm alg}.
$$
Set
$$
{\bP}_B^{\times_C k} =\prod^{k}_{i=1}{}_C {\bP}_B=
\underbrace{{\bP}_B\mathop{\times}\nolimits_C{\bP}_B
\mathop{\times}\nolimits_C\cdots\mathop{\times}\nolimits_C{\bP}_B}_{k}
$$
with the natural projection $\pi_B$ onto
$$
B^{\times_C k} = \prod^{k}_{i=1}{}_C B
=\underbrace{B\mathop{\times}\nolimits_CB\mathop{\times}\nolimits_C
\cdots\mathop{\times}\nolimits_C B}_{k},
$$
induced by the projection $\pi_B$ on each factor.

Let
$$
\tau\colon B^{\times_C k} \to
C^{\times_C k} =\prod^{k}_{i=1}{}_C C=C
$$
be the natural projection. Obviously, $\tau^{-1}(C^{\circ})$ is a
smooth affine algebraic variety. For an open (in any topology)
subset $D\subset C^{\circ}$ set
$$
Z_B(D)=\{(y_1,\dots,y_k)|\tau(y_1)=\dots=\tau(y_k)=x\in D,
y_i\neq y_j \,\,\, \mbox{for}\,\,\, i\neq j\}.
$$

Let
$$
Z_B\subset B^{\times_C k}
$$
be the closure of the open set $Z_B(C^{\circ})$. On the varieties
$B^{\mathop{\times}\nolimits_C k}$ and ${\bP}_B^{\mathop{\times}\nolimits_C k}$
we have the canonical action
of the symmetric group $S_k$: it permutes the factors of the
fiber product. This action is obviously compatible with the
projection $\pi_B$. For an open set $D\subset C^{\circ}$ let
$$
Z_{\bP}(D)=\pi^{-1}_B(Z_B(D))\subset {\bP}_B^{\mathop{\times}\nolimits_C k},
$$
and let
$$
Z_{\bP}\subset {\bP}_B^{\mathop{\times}\nolimits_C k}
$$
be the closure of the open set $Z_{\bP}(C^{\circ})$. The closed
sets $Z_B$, $Z_{\bP}$ are invariant with respect to $S_k$,
$\pi_B(Z_{\bP})=Z_B$, $\pi_B^{-1}(Z_B)=Z_{\bP}$, and the
projection $\pi_B$ is compatible with the action of the group
$S_k$. It is easy to see that the quotient variety satisfies
$$
B^{B/C}=Z_B/S_k\cong C.
$$
Now set
$$
{\bP}^{B/C}_{\rm alg}=Z_{\bP}/S_k.
$$

\begin{proposition} The algebraic variety ${\bP}^{B/C}_{\rm alg}$
and the compact complex space $\bP_{\rm an}^{B/C}$
are canonically isomorphic (as complex
spaces).
\end{proposition}

\proof  Let $D\subset C^{\circ}$ be a simple open set, so
that over $D$ the cover $\tau\colon B\to C$ is reducible,
$$
\tau^{-1}(D)=\coprod^k_{i=1}D_i.
$$
As above, the symbol $\pi$ denotes the natural projection for each
of the varieties
under consideration. By
definition,
$$
Z_{\bP}(D)=\{((p_1,y_1),(p_2,y_2),\dots,(p_k,y_k))|\tau(y_1)=\dots=\tau(y_k)=x\in
D, y_i\neq y_j\,\,\, \mbox{for}\,\,\, i\neq j\}.
$$
Consider the map
$$
\lambda\colon Z_{\bP}(D)\to{\bP}^{B/C}_{\rm an}(D)
$$
defined in the following way:
$$((p_1,y_1),\dots,(p_k,y_k))\mapsto
(p_{\sigma(1)};\dots;p_{\sigma(k)};z=\tau(y_1))\in{\bP}^{\times
k}\times D,
$$
where the permutation $\sigma\in S_k$ is uniquely determined by
the condition
$$
y_{\sigma(i)}\in D_i.
$$
It is easy to see that the fibers of the map $\lambda$ are
precisely the orbits of the group $S_k$, so that $\lambda$
defines a natural isomorphism
$$
Z_{\bP}(D)/S_k\cong {\bP}^{B/C}_{\rm an}(D).
$$
Since this isomorphism is natural, it extends to an isomorphism
$$
Z_{\bP}(C^{\circ})/S_k\cong {\bP }^{B/C}_{\rm an}(C^{\circ}).
$$
We must check that this isomorphism of complex spaces extends to
a global isomorphism
$$
Z_{\bP}/S_k\cong {\bP}^{B/C}_{\rm an}.
$$

To avoid overcomplicated notations, we perform this work for an
``extremal'' particular case of the maximal ramification, when $B$
has only one branch over $x\in C$. In the general case of several
branches the arguments are completely similar, only there are
more indices and the formulas are more complicated. So assume that
$\tau^{-1}(x)=\{y\}$ and with respect to suitable local
coordinates $z$ and $w$ on $C$ and $B$ the map $\tau$ is given by
the formula
$$
\tau^{-1}(D)\ni w\mapsto z=w^k\in D.
$$
Recall that
$$
\begin{array}{l}
{\bP}^{B/C}_{\rm an}(D)=\{(p_1;\dots;p_k;z)\}/{\zed}_k, \\ \\
{\bP}^{B/C}_{\rm alg}(D)=Z_{\bP}(D)/S_k,
\end{array}
$$
where $Z_{\bP}(D)$ is the closure of
$$
Z_{\bP}(D^{\circ})=\{((p_1,z_1),\dots,(p_k,z_k))|z_1^k=\dots=z_k^k=z\in
D^{\circ}, z_i\neq z_j \,\,\, \mbox{for}\,\,\, i\neq j\},
$$
$D^{\circ}=D\setminus \{x\}$ being a punctured disk.

Let the symbol $S^{\circ}({\zed}_k)$ denote the set (group) of
permutations of the additive group ${\zed}_k$, which maps 0 to
itself. We get the decomposition
$$
Z_{\bP}(D^{\circ})=\coprod_{\sigma\in S^{\circ}({\zed}_k)}Z_{\bP}(D^{\circ},\sigma),
$$
where
$$
Z_{\bP}(D^{\circ},\sigma)=\{((p_1,z_1),\dots,(p_k,z_k))|z_i=\zeta^{\sigma(i-1)}z_1\},
$$
$\zeta$ is a fixed primitive root of unit of degree $k$. The
group $S_k$ permutes the components $Z_{\bP}(D^{\circ},\sigma)$:
$$
\mu\colon Z_{\bP}(D^{\circ},\sigma)\to Z_{\bP}(D^{\circ},\mu[\sigma]),
$$
$\mu\in S_k$, where the permutation $\mu[\sigma]\in
S^{\circ}({\zed}_k)$ is defined by the relation
\begin{equation}
\label{4} \mu[\sigma](i-1)=\sigma(\mu(i)-1)-\sigma(\mu(1)-1),
\end{equation}
$i=1,\dots,k$, and we identify the integer $i$ and its image in
${\zed}_k$. Now set $Z_{\bP}(D,\sigma)$ to be the closure of
$Z_{\bP}(D^{\circ},\sigma)$ in $Z_{\bP}(D)$. Obviously, the
group $S_k$ permutes isomorphically the subsets $Z_{\bP}(D,\sigma)$
and acts transitively on the set of these subsets.
Therefore, if $N(\sigma)\subset S_k$ is the stabilizing subgroup
of the component $Z_{\bP}(D,\sigma)$, then
$$
{\bP}^{B/C}_{\rm alg}(D)\cong Z_{\bP}(D,\sigma)/N(\sigma)
$$
for each $\sigma\in S^{\circ}({\zed}_k)$. Take
$\sigma=\varepsilon$, the identical permutation. If $\mu\in
N(\varepsilon)$, then, according to (\ref{4}),
$$
\mu(i)-i \in{\zed}_k
$$
does not depend on $i=1,\dots,k$, that is, $\mu$ is a power of
the cycle
$$
(12\dots k).
$$
This fact gives us the required isomorphism
$$
{\bP}^{B/C}_{\rm alg}(D)\cong Z_{\bP}(D,\varepsilon)/{\zed}_k
\cong Z_{\bP}^{\rm an }(D)/{\zed}_k\cong {\bP}^{B/C}_{\rm
an}(D),
$$
since the coincidence $Z_{\bP}(D,\varepsilon)=Z_{\bP}^{\rm an
}(D)$ is obvious. It is easy to check that the isomorphism just
constructed is compatible with the isomorphisms
$$
{\bP}^{B/C}_{\rm alg}(D^{\sharp})\cong  {\bP}^{B/C}_{\rm
an}(D^{\sharp})
$$
that were constructed above for open subsets $D^{\sharp}\subset
D\cap C^{\circ}$.
\qed

\subsection{Divisors of fiber-wise type $(M,\dots,M)$ on ${\bP}^{B/C}$}

In this subsection we construct the special model $X^{\sharp}/C$ which was already discussed in
sect. 1.3.
Recall the rational map
$$
p\colon X\mathop{\times}\nolimits_C B-\,-\,\to {\bP},
$$
regular over $C^{\circ}$. Using the projection onto the second
factor, we extend $p$ to the map
$$
p_B\colon X\mathop{\times}\nolimits_C B-\,-\,\to {\bP}_B={\bP}\times B,
$$
also regular over $C^{\circ}$. Now construct the product
$$
p_B^{\mathop{\times}\nolimits_C k}\colon
\underbrace{(X\mathop{\times}\nolimits_C
B)\mathop{\times}\nolimits_C\dots\mathop{\times}\nolimits_C(X\mathop{\times}\nolimits_C
B)}_{k}-\,-\,\to\underbrace{{\bP}_B
\mathop{\times}\nolimits_C\dots\mathop{\times}\nolimits_C{\bP}_B}_{k}
$$
and note that there is a natural isomorphism
$$
(X\mathop{\times}\nolimits_C B)^{\mathop{\times}\nolimits_C
k}\cong X^{\mathop{\times}\nolimits_C
k}\mathop{\times}\nolimits_C B^{\mathop{\times}\nolimits_C k}.
$$
Taking the composition with the diagonal embedding
$$
\Delta\colon
X\hookrightarrow X^{\mathop{\times}\nolimits_C k},
$$
we get the rational map
$$
p_B^{\mathop{\times}\nolimits_C k}\circ
(\Delta,\mathop{\rm id})\colon X\mathop{\times}\nolimits_C
B^{\mathop{\times}\nolimits_C k}-\,-\,\to {\bP}_B^{\mathop{\times}\nolimits_C k}.
$$
It is easy to see that this map is equivariant with respect to
the action of the symmetric group $S_k$ and that the image of
$X\mathop{\times}\nolimits_C Z_B$ lies in $Z_{\bP}$.
Factorizing by $S_k$ and taking into consideration that $Z_B/S_k$
is just $C$, we get a map
$$
\iota\colon X-\,-\,\to {\bP}^{B/C},
$$
where $\iota|_{X(D)}=\iota_D\colon X(D)\to{\bP}^{B/C}(D)$ for
any $D\subset C^{\circ}$. The image
$$\iota(X)=X^{\sharp}\subset{\bP}^{B/C}$$ is the desired special
model of $X/C$.

Thus we have realized (on the birational level) the fiber space
$X/C$ as a hypersurface $X^{\sharp}\subset{\bP}^{B/C}$ of
fiber-wise type $(M,\dots,M)$. Now let us give an explicit
description of such hypersurfaces.

Fix a system of homogeneous coordinates $(u_0,\dots,u_M)$ on
${\bP}$. Let
$$
{\bf M}=\{{\cal M}=u_0^{a_0}u_1^{a_1}\dots
u_M^{a_M}|a_0+\dots+a_m=M\}=\{{\cal M}_j|j\in J\}
$$
be the set of all monomials of degree $M$ in the variables $u_*$.
For every point $x\in C^{\circ}$ the fiber
$\tau^{-1}(x)=\{y_1,\dots,y_k\}$ consists of precisely $k$
points, and the fiber ${\bP}^{B/C}_x$ naturally identifies with
the product ${\bP}^{(1)}\times\dots\times {\bP}^{(k)}$, where
${\bP}^{(i)}={\bP}$ corresponds to the point $y_i\in B$.
Consider the set of all monomials of the type $(M,\dots,M)$ on the
space ${\bP}^{B/C}_x$:
$$
{\bf M}^{\tau^{-1}(x)}=\{{\cal M}^{(1)}_{i_1}\dots {\cal
M}^{(k)}_{i_k}|i_l\in J\},
$$
where the monomial ${\cal M}^{(l)}_{i_l}$ is in the variables
$u^{(l)}_j$, the homogeneous coordinates on the space ${\bP}^{(l)}$. Obviously, there are
$$
\sharp ({\bf M}^{\tau^{-1}(x)})=(\sharp {\bf M})^k
$$
such monomials.

Defining a hypersurface of the type $(M,\dots,M)$ in ${\bP}^{\times k}={\bP}^{B/C}_x$
is equivalent to give a point in
the projective space
$$
{\bP}^{(\sharp M)^k-1}
$$
or a non-zero set of coefficients
$$
\{c_{i_1,\dots,i_k}\in{\comp}|i_l\in J\}/{\comp}^*
$$
modulo a non-zero factor. For a given hypersurface of the type
$(M,\dots,M)$ we can change the system of coordinates and make
sure that for some index $e\in J$
$$
c_{e,e,\dots,e}\neq 0.
$$
Normalizing this coefficient to 1, we determine all the
other coefficients $c_{i_1,\dots,i_k}$.

In particular, for our hypersurface $X^{\sharp}\subset{\bP}^{B/C}$ we obtain:  \\
at a point $x\in C^{\circ}$ of general position
the set of coefficients, defining the hypersurface $X^{\sharp}_x$
in ${\bP}^{B/C}_x$, is unique. It is easy to see that the
functions on $C^{\circ}$ that we get as coefficients are branches
of multi-valued rational algebraic functions on $C$. More
precisely, let $D\subset C^{\circ}$ be a simple open set
$$
\tau^{-1}(D)=\coprod^k_{i=1}D_i\subset B.
$$
Then we get a unique set of meromorphic functions
$$
c_{j_1,\dots,j_k}\colon D\to{\comp},
$$
defining the hypersurface $X(D)$:
$$
X(D)=\left\{ \sum_{(j_1,\dots,j_k)\in J^k} c_{j_1,\dots,j_k} {\cal
M}^{(1)}_{j_1}\dots {\cal M}^{(k)}_{j_k}=0\right\}.
$$
Now we have to globalize the meromorphic functions $c_{\sharp}$
as branches of multi-valued rational functions on $C$, that is,
rational functions on certain covers of the curve $C$. Obviously,
if
$$
c_{i_1,\dots,i_k}\quad \mbox{and}\quad c_{j_1,\dots,j_k}
$$
are two branches of one and the same multi-valued meromorphic
function, then the $k$-uples
$$
(i_1,\dots,i_k)\quad \mbox{and}\quad (j_1,\dots,j_k)
$$
differ by a permutation of the symbols only. More formally, for
any partition
$$
k=k_1+\dots+k_r
$$
and any set
$$
{\cal M}_1,\dots,{\cal M}_r
$$
of pair-wise distinct monomials from ${\bf M}$ consider all the
maps
$$
\varepsilon\colon\tau^{-1}(x)\to{\bf M},\quad
\sharp\varepsilon^{-1}({\cal M}_i)=k_i.
$$

Such maps and only such maps (for a fixed partition and a fixed
set of pair-wise distinct monomials ${\cal M}_i$) can generate
different branches of the same multi-valued function.

Let
$$
\alpha_{k_1,\dots,k_r}\colon A_{k_1,\dots,k_r}\to C
$$
be the finite cover, parametrizing for any point $x\in C^{\circ}$
of general position all the partitions of the fiber $\tau^{-1}(x)$
into a system of subsets of cardinality $k_1,\dots,k_r$:
$$
A_{k_1,\dots,k_r}=Z_B/(S_{k_1}\times\dots\times S_{k_r}),
$$
where $S_{k_1}\times\dots\times S_{k_r}\hookrightarrow S_k$ is
the natural inclusion. The curve $A_{k_1,\dots,k_r}$ certainly
can be reducible. In this case when we speak about a rational
function
$$
a_{k_1,\dots,k_r}\colon A_{k_1,\dots,k_r} \to{\comp}
$$
we mean an arbitrary set of rational functions on its components.

A hypersurface of fiber-wise type $(M,\dots,M)$ in ${\bP}^{B/C}$
is defined by a set of rational functions
$$
a^{{\cal M}_1,\dots,{\cal M}_r}_{k_1,\dots,k_r}\colon
A_{k_1,\dots,k_r}\to{\comp}
$$
for all the partitions $k=k_1+\dots+k_r$, $k_1\geq \dots\geq
k_r\geq 1$, and all the $r$-uples of pair-wise distinct monomials
${\cal M}_1,\dots,{\cal M}_r$ from ${\bf M}$. The sets
$$
\left\{a^{{\cal M}_1,\dots,{\cal
M}_r}_{k_1,\dots,k_r}\right\}\quad\mbox{and}\quad \left\{b^{{\cal
M}_1,\dots,{\cal M}_r}_{k_1,\dots,k_r}\right\}
$$
define the same hypersurface if and only if there exists a
non-zero rational function
$$
\lambda\colon C\to{\comp}
$$
such that
$$
a^{{\cal M}_1,\dots,{\cal M}_r}_{k_1,\dots,k_r}=(\lambda\circ
\alpha_{k_1,\dots,k_r})b^{{\cal M}_1,\dots,{\cal
M}_r}_{k_1,\dots,k_r}
$$
for all $(k_1,\dots,k_r)$, $({\cal M}_1,\dots,{\cal M}_r)$. \\

Fixing a set of rational functions $a^{\sharp}_{\sharp}$, which defines
the hypersurface $X^{\sharp}=\iota(X)\subset{\bP}^{B/C}$,
our description of the special model $X^{\sharp}\subset{\bP}^{B/C}$
is complete: everything is ready for investigation of
the problem of existence of a section of the morphism $\pi\colon
X^{\sharp}\to C$.

\section{Existence of sections}

In accordance with our strategy, we look at the sections of the
fibration $X^{\sharp}/C$ as the sections of the ambient fibration
$s\colon C\to{\bP}^{B/C}$, $\pi\circ s=\mathop{\rm id}_C$, the
image of which lies in $X^{\sharp}$, i.e. $s(C)=X^{\sharp}$.

\subsection{Sections of the ambient fibration ${\bP}^{B/C}$}

Let $x\in C^{\circ}$. To give a section $s\colon C\to{\bP}^{B/C}$ means to give a value
$$
s(x)\in {\bP}^{B/C}_x=\prod_{b\in\tau^{-1}(x)}{\bP}_b
$$
for every point $x\in C^{\circ}$, which is the same as to give
values
$$
\{s_B(b)\in{\bP}_b|b\in\tau^{-1}(x)\}.
$$
Hence sections of the ambient fibration are in one-to-one
correspondence with sections
$$
s_B\colon B\to {\bP}_B={\bP}\times B.
$$
In their turn, the latter sections are determined by $(M+1)$-uples
of rational functions on $B$.

Now we introduce certain finite-dimensional spaces of sections in
the following way. Let $R$ be an effective divisor on $B$. Set
$$
{\cal S}(R)={\cal L}(R)^{M+1}=\{(s_0,\dots,s_M)|s_i\in{\cal
L}(R) = H^0({\cal O}_B(R))\}.
$$
Obviously,
$$
\dim{\cal S}(R)=(M+1)\dim{\cal L}(R)
$$
and any non-zero $s\in{\cal S}(R)\setminus\{0\}$ determines a
section of the fibration ${\bP}_B/B$ and thus of ${\bP}^{B/C}/C$. If $\mathop{\rm deg}R\gg 0$ is, then the map
$$
{\bP}({\cal S}(R))\to \mathop{\rm Sect}({\bP}_B/B)=\mathop{\rm
Sect}({\bP}^{B/C}/C)
$$
is injective on a non-empty Zariski open subset of ${\bP}({\cal
S}(R))$.

>From now on we assume that $\tau(R)\subset C^{\circ}$, so that
$R$ does not contain ramification points of the morphism $\tau$.

\subsection{Sections of the fibration $X^{\sharp}/C$}

Let $x\in C^{\circ}$, $D\ni x$ be a simple open set, so that
$$
\tau^{-1}(D)=\coprod^k_{i=1}D_i,
$$
$\tau\colon D_i\to D$ are isomorphisms for all $i$. Consider
$$
{\bP}^{B/C}(D)={\bP}^{(1)}\times\dots\times{\bP}^{(k)}\times D,
$$
and let $u^{(i)}_j$ be the homogeneous coordinates on ${\bP}^{(i)}$,
which come from the coordinates $u_{\sharp}$ on ${\bP}$. In terms of these coordinates the hypersurface
$X^{\sharp}(D)\subset{\bP}^{B/C}(D)$ is given by the equation
$$
F_D(u^{(1)}_i;\dots;u^{(k)}_i;z)=\sum_{(i_1,\dots,i_k)\in J^k}
c_{i_1,\dots,i_k}(z){\cal M}^{(1)}_{i_1}\dots{\cal M}^{(k)}_{i_k},
$$
where the monomials ${\cal M}^{(j)}_i$ are in the variables
$u^{(j)}_{\sharp}$ and of degree $M$, whereas the meromorphic
functions $c_{\sharp}(z)$ are branches of the functions
$$
a^{{\cal M}_1,\dots,{\cal M}_r}_{k_1,\dots,k_r},
$$ with
$z$ a coordinate on $D$. The single-valued branches exist
because all the maps
$$
\alpha_{k_1,\dots,k_r}\colon A_{k_1,\dots,k_r}\to C
$$
are reducible over $D$:
$$
\alpha_{k_1,\dots,k_r}^{-1}(D)=\coprod_{\beta\in{\cal
B}(k_1,\dots,k_r)}D_{\beta},
$$
where ${\cal B}(k_1,\dots,k_r)$ is the set of all partitions of
the type $(k_1,\dots,k_r)$ of the set $\{1,\dots,k\}$,
parametrizing the points of the fiber $\tau^{-1}(x)$.

Take $s=(s_0,\dots,s_M)\in{\cal S}(R)$. Obviously,
$$
s|_{\tau^{-1}(D)}=\{s|_{D_i}=(s^{(i)}_0,\dots,s^{(i)}_M)|i=1,\dots,k\}.
$$
The image of the section $s_D\colon D\to {\bP}^{B/C}(D)$,
defined by $s\in{\cal S}(R)\setminus\{0\}$, lies in
$X^{\sharp}(D)$ if and only if
$$
F_D(s^{(1)}_{\sharp};s^{(2)}_{\sharp};\dots;s^{(k)}_{\sharp};z)\equiv
0.
$$
Let us understand what happens when $s^{(i)}_{\sharp}$ are
substituted as $u^{(i)}_{\sharp}$ in the equation $F$. Let
$s\in{\cal S}(R)$ be an arbitrary element (possibly zero).

\begin{proposition} The meromorphic functions of the
variable $z$
$$
F_D(s_D)=F_D(s^{(1)}_{\sharp};\dots;s^{(k)}_{\sharp};z)
$$
glue together into a single-valued rational algebraic function
$F(s)$ on $C$.
\end{proposition}

\proof  This is obvious, since both $s^{(i)}_j(z)$ and
$c_{\sharp}(z)$ come from the globally (over $C$) defined objects:
the branches of multi-valued algebraic functions on $C$, and for
$D_1\subset D_2\subset C^{\circ}$ we get
$F_{D_1}(s_{D_1})=F_{D_2}(s_{D_2})|_{D_1}$, so that $F_D(s_D)$
glue together into a single-valued function on $C^{\circ}$ and
thus on $C$.
\qed

\subsection{The crucial fact}

Let $(F(s))$ be the divisor of poles and zeros of the rational
function $F(s)$ on $C$.

\begin{proposition} There exists an effective divisor
$\Delta$ on the curve $C$, depending on the hypersurface
$X^{\sharp}$ only (but not depending on $R$), such that for
$s\in{\cal S}(R)$
\begin{equation}
\label{5} (F(s))+M\tau_*R+\Delta\geq 0.
\end{equation}
\end{proposition}

\proof We observe that up to a divisor $\Delta$
which does not depend on $R$ the poles of the rational function
$F(s)$ are bounded from above by the divisor
$$
M\tau_*R=M\left(\sum_{i\in I} m_i\tau(y_i)\right),\quad
\mbox{where}\quad R=\sum_{i\in I} m_i y_i.
$$
Let us prove first a local version of the estimate (\ref{5}). Let
$D\subset  C^{\circ}$ be a simple open set,
$$
\tau^{-1}(D)=\coprod^k_{i=1}D_i,
$$
$\tau\colon D_i\to D$ is an isomorphism, $z$ is a coordinate on
$D$, $z_i=\tau^* z|_{D_i}$ is a coordinate on $D_i$. Consider the
meromorphic function
$$
F_D(s_D)=F_D(s^{(1)}_{\sharp};\dots;s^{(k)}_{\sharp};z)=
$$
$$
=\sum_{(i_1,\dots,i_k)\in J^k} c_{i_1,\dots,i_k}(z){\cal
M}_{i_1}(s^{(1)}_{\sharp}\dots{\cal M}_{i_k}(s^{(k)}_{\sharp}),
$$
where $s^{(i)}_j(z)=(\tau^{-1})^*(s_j|_{D_i})$.\\
The local version of (1) is now provided by

\begin{proposition} There exists an effective divisor
$\Delta^+$ on $C^{\circ}$, depending on $X^{\sharp}$ only, such
that for any simple open set $D\subset C^{\circ}$
$$
(F_D(s_D))+M\tau_*(R\cap \tau^{-1}(D))+\Delta^+\cap D\geq 0.
$$
\end{proposition}

\proof  Since the degree of the monomials ${\cal M}_i$ is
equal to $M$, whereas by construction for any $i,j$
$$
(s^{(i)}_j)+(R\cap D_i)\geq 0,
$$
we get for any $j\in J$
$$
({\cal M}_j(s^{(i)}_{\sharp}))+M(R\cap D_i)\geq 0.
$$
This implies
$$
\left(\prod^k_{e=1}{\cal
M}_{i_e}(s^{(e)}_{\sharp})\right)+M\tau_*(R\cap
\left(\coprod^k_{e=1}D_i\right))\geq 0,
$$
that is,
$$
\left(\prod^k_{e=1}{\cal
M}_{i_e}(s^{(e)}_{\sharp})\right)+M\tau_*(R\cap \tau^{-1}(D))\geq
0.
$$
Recall that the functions $c_{i_1,\dots,i_k}$ are branches of the
rational functions
$$
a^{{\cal M}_1,\dots,{\cal M}_r}_{k_1,\dots,k_r}\colon
A_{k_1,\dots,k_r}\to{\comp}
$$
and thus depend on $X^{\sharp}$ only. Their poles are fixed and
therefore there exists an effective divisor $\Delta^+$ on
$C^{\circ}$ such that
$$
(F_D(s_D))+M\tau_*(R\cap \tau^{-1}(D))+\Delta^+\cap D\geq 0,
$$
which is what we need. \qed

 Now we globalize (3.3) and prove Proposition 3.2
From the local fact which we have just obtained we immediately
get the estimate
$$
(F(s))\cap C^{\circ}+M\tau_*R+\Delta^+\geq 0
$$
(since $\mathop{\rm Supp}\tau_*R\subset C^{\circ}$ by assumption).
Consequently, to complete the proof of Proposition 3.2 it remains
to consider the poles of the function $F(s)$ at the points $x\in
C\setminus C^{\circ}$ and prove that they are bounded. But this
is easy and done as follows.

In order to be able to work with coordinates, we use the analytic
description of ${\bP}^{B/C}$ near a point $x\in C\setminus
C^{\circ}$ given above. Let $D\ni x$ be a small disk, then
$$
{\bP}^{B/C}(D)=Z_{\bP}(D)/G_1\times\dots\times G_a.
$$
Let us pull back all the functions to $Z_{\bP}^{\rm an}(D)$.
The polynomial $F_D$ outside the point $x$, that is, on
$D^{\sharp}=D\setminus\{x\}$, is defined by the construction,
given above. Since by assumption $R\subset C^{\circ}$, all the
functions $s^{(i)}_j$ are regular on $\tau^{-1}(D)$, and therefore
the monomials ${\cal M}_j$ in these functions are regular, too.
The polynomial $F_D$ is linear in its coefficients, which depend
on $X^{\sharp}$ only. Thus there exists an effective divisor
$\Delta^{\sharp}$ with the support in $C\setminus C^{\circ}$,
$$
\Delta^{\sharp}=\sum_{x\in C\setminus C^{\circ}}m_xx,
$$
such that near $x\in C\setminus C^{\circ}$ we have
$$
(F_D(s))+(\Delta^{\sharp}\cap D)\geq 0.
$$
Setting $\Delta=\Delta^++\Delta^{\sharp}$, we complete the proof
of Proposition 3.2.
\qed

\subsection{End of the proof of Theorem 1.2}

As we have seen above, the map
$$
\varphi\colon s\mapsto F(s)\in {\comp}(C)
$$
is a map between affine varieties
$$
\varphi\colon{\cal S}(R)\to{\cal L}(M\tau_*R+\Delta).
$$
By construction (locally over each simple open subset $D\subset
C^{\circ}$ and therefore globally due to the well-known
characterization of polynomial functions) $\varphi$ is polynomial
of degree $Mk$, that is, makes a regular map of affine spaces
$$
\varphi\colon{\comp}^{\dim{\cal S}(R)}\to{\comp}^{\dim{\cal
L}(M\tau_*R+\Delta)}.
$$
Obviously $0\in\varphi^{-1}(0)$, so that the fiber
$\varphi^{-1}(0)$ over the point $0\in{\cal L}(M\tau_*R+\Delta)$
is non-empty. Consequently,
\begin{equation}
\label{6}
\begin{array}{ccc}
\dim\varphi^{-1}(0) & \geq & \dim{\cal S}(R)-\dim{\cal L}(M\tau_*R+\Delta) \\
                          &   & \|   \\
                          &    &  \mathop{\rm deg}R-K,
\end{array}
\end{equation}
where $K\in{\zed}$ is a constant. Therefore if $\mathop{\rm
deg}R\gg 0$ is sufficiently large, then the fiber
$\varphi^{-1}(0)$ contains non-zero elements, which define
sections of the fibration $X^{\sharp}/C$.\\

The proof of Theorem 1.2 is now complete. \\

The general theory [KMM96] now implies that the fibration
$X^{\sharp}/C$ has plenty of sections, connecting arbitrary points
$w_i\in X^{\sharp}$ in distinct fibers. The map $\varphi$
constructed above gives an explicit description of sections as
solutions of a system of algebraic equations, whereas the
inequality (\ref{6}) gives an estimate for the dimension of the
space of sections with bounded poles. The concluding part of our
arguments is a generalization of the classical proof of the Tsen
theorem for the ``twisted'' (not locally trivial) fibration ${\bP}^{B/C}\to C$.

\section{Fibrations of Fano threefolds}

Our first result is

\begin{theorem} Let $\pi: X \to C = \bP_1$ be a fibrations with general fiber $X_b$ a Fano $3-$fold with $b_2(X_b) = 1.$ Then
$\pi$ has a section.
\end{theorem}

\proof We are going to apply the classification of Fano $3-$folds and refer to [IP99] with the further references given there.
Let $r$ be the index of $X_b.$ \\
If $r = 4$ or $r =  3,$ then $X_c $ is $\bP_3$ or $Q_3,$ and the claim follows from Theorem 1.2. \\
If $ r = 2,$ then $X_b$ is a del Pezzo $3-$fold and all cases are again settled by 1.2 except for $X_b$ being of type $V_5,$ i.e.
$X_b$ is a linear section of $G(2,5)$ embedded by Pl\"ucker in $\bP_9.$ In that case we consider the relative Douady space
$\tau: \sD \to B$ of lines in $\pi-$fibers. Since the surface of lines in $X_c$ is just $\bP_2,$ the map $\tau$ is a generic $\bP_2-$
fibration over $B.$ Thus $\tau $ has a section providing a single rational curve in the general fiber $X_c.$ These rational
curves define a rational surface in $X$ over $\bP_1,$ so that $\pi$ has a rational section. \\
It remains to treat the case $r = 1.$ Let $g$ be the genus of $X_c.$ If $g \leq 5,$ then again $X_c$ is a complete intersection in a
weighted projective space, so we can apply 1.2. So let $g \geq 6;$ here we will use the Mukai construction. Let $C^* \subset C$ be the set
of regular values of $\pi$ and fix $c  \in C.$ Then there exists a unique stable vector bundle $\sE_c$ over $X_c$ which is spanned
and whose sections define a map
$$ \phi_c: X_c \to G(m,V_c)$$
where $V_c = H^0(\sE_c)$ such that $\sE_c = \phi^*(E)$ for the tautological bundle $E$ on $G(m,V_c).$ Note that $\phi_c$ is
an embedding if $g \geq 7$ but it may be a two-sheeted cover if $g = 6.$ Moreover $\phi_b(X_c)$ is a complete intersection
in $G(m,V_c)$ if $g \geq 8.$ Now every stable bundle on $X_c$ of the same rank and same Chern classes of $\sE_c$ is isomorphic
to $\sE_c.$ The $\sE_c$ clearly glue locally; and we are going to show that they even glue globally to a vector bundle $\sE$ over
$X^* = \pi^{-1}(C^*).$
In fact, choose an open covering $U_i \subset  C^*$ such that $\sE_i$ exists over $V_i = \pi^{-1}(U_i).$ Then the uniqueness of $\sE_c$
yields isomorphisms
$$ \phi_{ij}: \sE_i \vert V_{ij} \to \sE_j \vert V_{ij}.$$
By stability, we have ${\rm Aut}(\sE_c) = \comp^*,$ hence the cocycle condition
$$ \phi_{ij} \phi_{jk} = \phi_{ik} $$
is satisfied up to $\comp^*.$ This means that at least that the $\bP(\sE_i)$ glue to a projective bundle over $X^*.$ Since $X$ is a rationally
connected fibration over $\bP_1,$
this bundle has to come from a vector bundle $\sE.$ \\
As a consequence we obtain a map
$$\psi: X \rightharpoonup Z $$
to a fibration $\sigma: Z \to C$ which is a Grassmann bundle over $C^*$ and which is birationally a product.
Moreover $\psi \vert X_c = \phi_c$ for $c  \in C^*.$
If now $g = 8 $ or $g = 6$ and $X$ is not special, then we can apply the next proposition (Tsen for Grassmann bundles) to
conclude. The special case (with genus 6) is the following: the sections of $\sE_c$ provide a $2-$sheeted covering
$$ \phi_c: X_c \to G(2,5) $$
the image of which is a Fano $3-$fold of type $V_5.$ Now the (analytic) preimage of a line in $V_5$ is a conic
in $X_c$ and conversely every conic in $X_c$ is mapped to a line in $V_5.$ Therefore the conics in $X_c$ form an irreducible
family, parametrized by a surface $S$, which admits a map $h: S \to \bP_2,$ the surface of line in $V_5.$ Clearly
two smooth conics cannot be mapped to the same line, hence $h$ is birational, i.e. $S$ is rational. Now the same
argument as in the case $V_5$ provides a section in $X.$ \\
In the remaining cases more efforts are necessary, given in 4.3 and 4.5 below.
\qed

\begin{proposition} Let $G = G(p,N)$ be a Grassmannian (projective notation), let $B$ be a smooth rational curve or even smooth curve of any genus $g$
and let $H \subset B \times G$ are relative smooth complete intersection of elements
$$  H_i \in \vert \sO_G(k_i) \otimes g^*(\sL_i)\vert. $$
If $\sum k_i$ is smaller as the index of $G$, i.e. $\sum k_i < N+1,$ then $H \to B$ has a section.
\end{proposition}

\proof
For simplicity of notations we treat only the case where $H$ is a hypersurface and leave the trivial amplification
to the general case to the reader. \\
A section $s$ of $H \to B$ is given by $p+1$ elements $s_i \in H^0(B,(\sO_B(d_i))^{\oplus(N+1)}),$ where $0 \leq i \leq N+1,$
subject to the following conditions:
\begin{enumerate}
\item $s_0(b), \ldots , s_p(b) $ are linearly independent for general $b \in B$;
\item  $s_0(b) \wedge \ldots \wedge s_p(b) \in H_b$ for all $b  \in B.$
\end{enumerate}

Set
$$
V(d_*) := \mathop{\bigoplus}\limits^p_{i=0}
H^0(\sO_B(d_i)^{\oplus (N+1)})
$$
to be the space parametrizing the $(p+1)$-uples of elements $s_i$
and consider the map
$$
\pi: V(d_*) \to \mathop{\bigoplus}\limits_{
\begin{array}{c}
I\subset \{1,\dots,N\} \\
|I|=p+1
\end{array}
} H^0(\sO_B(d)) =: P_d
$$
which is given by Pl\" ucker coordinates and where we set
$$
d=\sum^p_{i=0}d_i.
$$
Its image may well be non-closed. Denote by $W(d_*)$ its closure
in $P_d$. The first problem is to compute its dimension $\dim W(d_*)$.
From now on we will assume that $d_0 = d_1 = \ldots = d_{p} = e$,
so that $d=(p+1)e$. It is enough to prove our claim for one particular
choice of degrees $d_i$, so let us take this one. In this case denote
$V(d_*)$ by $V_e$, $W(d_*)$ by $W_e$. Note that the morphism
$$
\pi\colon V_e\to W_e
$$
is dominant by construction.
\\

{\bf Lemma.} {\it The relative dimension
$$\dim V_e/W_e = \dim V_e - \dim W_e =  (p+1)^2
$$
as soon as $e \geq 2g+1$.}

Putting off the proof of the lemma for a while, let us complete the
proof of (4.2). Indeed, if the hypersurface $H$ is defined
by the equation $F \in H^0(\sO_G(k) \otimes \pi^*(\sL)),$
then the condition (2) is rewritten in Pl\"ucker coordinates as
$$
F(s(b)) = \sum_{ \vert \alpha \vert = k} p_{\alpha}(b) s^{\alpha} = 0.
$$
A section $s=(s_I)$ of the trivial fibration $B\times G/B$ given by its
Pl\" ucker coordinates $s_I$ is actually a section of the fibration $H/B$ if and
only if $F(s(b))\equiv 0$. Now assume that $s\in W_d$. so that
$s_I\in H^0(\sO_B(d))$ for each $I\subset \{1,\dots,N\}$. It is clear
that the map
$$
s\mapsto F(s)\in H^0(\sO_B(kd+m))
$$
where $m = {\rm deg}(\sL)$, determines a polynomial map
$$
\varphi\colon W_d\to H^0(\sO_B(kd+m))
$$
of closed affine algebraic varieties. Now $\varphi^{-1}(0)\ni 0$,
so that the fiber over $0$ is non-empty. Thus
$$
\dim \varphi^{-1}(0)\geq \dim W_d-h^0(\sO_B(m+kd))\sim
$$
$$
\sim \dim V_d-kd\sim (N+1-k)d
$$
up to a bounded integer. Thus, precisely as in the usual
Tsen theorem, we obtain a section. \qed for (4.2).
\\

{\bf Proof of the lemma.} We have to show that if $(s_0, \ldots, s_p)$ is as above and generic and if
$(s'_0, \ldots, s'_p),$ satisfies
$$ s'_0 \wedge \ldots \wedge s'_p = s_0 \wedge \ldots \wedge s_p;$$
so that
$$ s' = As \eqno (+) $$
with a $(p+1)\times (p+1)-$matrix $A = (a_{ij})$ of meromorphic functions on $B$ with $\det A = 1$, then $A$ is constant,
i.e. all $a_{ij} $ are constant functions.
To this extent we define a meromorphic map
$$ \psi: B \times \bP_p \rightharpoonup \bP_N$$
by setting $$\psi(b,[a_0:a_1: \ldots :a_p]) = [a_0s_0+ \ldots a_ps_p(b)] = [\sigma(b)].$$
Here the bracket $[\ . \ ]$ on the right hand side means the projection from a vector space $V \setminus 0$
to projective space, resp. $[\sigma(b)] = [\sigma_0(b): \ldots :\sigma_N(b)],$ where the $\sigma_j$
are the components of $\sigma = \sum a_i s_i.$
Now we will show
\\

{\bf Claim (4.2.2)} If $e\geq 2g+1$ is sufficiently large, then
$\psi $ is holomorphic and birational onto its image.
\vspace{0.3cm}

\noindent
First we show how to deduce our lemma from (4.2.2).
For simplicity of notation, let $\sG = \sO_B(e).$ Consider again
$s'\in H^0(B,\sG^{\oplus N+1})$ such that $s' = As$ with $\det A = 1.$ We shall write
$s' = \sum a_j s_j$.
If $e \geq 2g+1,$ then $\sG$ is very ample and therefore
the map $\Phi_{s''}: B \to \bP_N$ mapping $b$ to $[s''(b)],$ is an embedding for general $s'' \in  H^0(B,\sG^{\oplus N+1}) $ and the degree of the curve
$\Phi_{s''}(B) \subset \bP_N$ is independent of $s''.$ Since $s$ is generic and since we only have to
compute the dimension of the fiber at the general point, we may assume that $s'$ is general in the sense
that the previous sentences apply to $s' $ instead of $s''.$
Since $s' = As,$ we conclude that
$\Phi_{s'}(B) \subset W,$ the image of the birational map $\psi.$ So if $B'$  denotes the image of the
map $$ B \to B \times \bP_p, \ b \mapsto (b,[s'(b)]),$$
then $\Phi_{s'}(B) = \psi(B').$  Now in homology we can write
$$ [B'] = [B \times x] + a [l] $$
with a line $l \subset \bP_p$ and some integer $a \geq 0.$ Let $H = \psi^*(\sO_{\bP_N}(1)).$
Then $H \cdot B' = H \cdot B + a.$ On the other hand, $H \cdot B' = \deg \Phi_{s'}(B) =
\deg \Phi_s(B) = H \cdot B.$ Therefore $a = 0$ and there exists
$a = [a'_0: \ldots a'_p] \in \bP_p$ such that
$ [s'] = [\sum a'_j s_j]$. Hence there is a meromorphic function $a$ on $B$ such that
$$ s' = a (\sum a'_j s_j),$$
i.e. $ a a'_j = a_j $ for all $j.$
We want to argue that then $a$ must be constant. If not, then $s'$ would always have zeroes. But
$s'$ is a general section of a direct sum of very ample line bundle over a curve and therefore
has no section.
\\
\\ {\it Proof of (4.2.2)} If $d' \geq 2g+1,$
then the evaluation map
$$ \epsilon: H^0(\sG^{\oplus N+1}) \to \sG_b^{\oplus N+1} \oplus \sG_{b'}^{\oplus N+1} $$
is surjective for $b,b' \in B.$ This includes as usually the case $b = b',$ in which case the sum is understood as
$ \sG_b^{\oplus N+1} \otimes \sO_{B,b}/{m}_b^2.$
Therefore we obtain an embedding $$\tilde \psi: B \times \bP_p \to \bP_{N'}$$
with $N' = h^0(B,\sG^{\oplus N+1}) - 1.$ Recall that - possibly after dualizing the Grassmannian - we may assume that
$$2(p+1) \leq N+1.$$
Now let $\tau: \bP_{N'} \to \bP_N$ be the projection from a generic $\bP_M$ with $M = N'-N-1.$ Then
$\tau $ is holomorphic since $N+1 > p+1.$ Moreover $\psi$ is of the form $\psi = \tau \circ \tilde \psi$ and $\psi$ is birational unless
every $\bP_{N'-N}$ meets $\tilde \psi (B \times \bP_p)$ at least twice.  This can happen only if
$$p+1 = \dim B \times \bP_p
\geq {\rm codim}_{\bP_N'} \bP_{N'-N}, $$
i.e. if $p+1 \geq N.$ Since $N \geq 2(p+1)-1,$ as stated above, it follows $p = 0$ and
$N = 1,$ a trivial case. This finishes the proof of (4.2).
\qed

We are now attacking the cases $g = 9,10,12$. The key is (4.3); here it is more convenient to use affine notations
for the Grassmannians.

\begin{proposition} Let $G = G(p,m) $be the Grassmannian of $p-$dimensional linear subspaces in $\comp ^m.$
Let $i: G \to \bP^{{m \choose p}-1}$ be the Pl\"ucker
embedding and $B$ be a smooth rational curve. Let $F \subset B \times \bP^{{m \choose p}-1}$ be a relative linear subspace. Let $\sigma_1, \ldots,
\sigma_k$ be relative skew-symmetric $l-$linear forms on $B \times \comp^m$ satisfying simultaneously one of the following conditions
\begin{enumerate}
\item $l = 2$;
\item $ l = p-1$ and $k = 1.$
\end{enumerate}
Let $\Sigma \subset B \times G$ be the locus of $p-$planes which are $\sigma_j-$isotropic for $1 \leq j \leq k.$
Let $Z = F \cap \Sigma;$ we assume that $Z$ is of the expected dimension and smooth over the general point of $B.$ Furthermore let $n = \dim Z - 1.$
Assume finally that in case (1) we have $$m+n > p(m-p) - {{(p-1)} \choose 2} k $$ and that in case (2) we have $$m+n > p(m-p)-1.$$
Then the projection $Z \to B$ has a section.
\end{proposition}

{\bf (4.4)} We first show how to derive Theorem 4.1 in cases $g = 9,10,12$ from Proposition 4.3.  The first part of the proof of Theorem 4.1
gives a bimeromorphic map
$\psi : X \rightharpoonup Z $ with $Z \subset C \times G \subset C \times \bP^{{m \choose p}-1}.$ Here the relevant data are:
\begin{itemize}
\item $(m,p) = (6,3)$ if $g = 9$;
\item $(m,p) = (7,5)$ if $g = 10,$
\item $(m,p) = (7,3) $ if $g = 12.$
\end{itemize}
Moreover from Mukai's classification, see
[IP99,5.2.3], $Z$ is of the form $Z = F \cap \Sigma$ for appropriate $F$ and $\Sigma$ with corresponding data (in the notation of 4.3):
\begin{itemize}
\item $(l,k) = (2,1)$ if $g = 9,$
\item $(l,k) = (4,1)$ if $g = 10 $
\item $(l,k) = (2,3)$ if $g = 12.$
\end{itemize}
Now it is just a matter of direct calculation to verify the inequalities we need to apply 4.3. \\

\noindent {\bf Proof} {\it of Proposition 4.3} \
The section $s$ we are looking for will be given by a $p-$tuple of sections $(s^j)_{1 \leq j \leq p}: C \to \comp ^m,$ given by relative
polynomials $s^j$ on $C$ of degree $d$ - at the end we will take limits as $d \to \infty.$
So our section $s$ depends on $mp(d+1)$ parameters. \\
{\it Case (1)} ($l = 2)$:  The section $s$ we are looking for is subject to the following conditions:
$$ \sigma_i(s^j,s^l) = 0 $$
for $1 \leq i \leq k$ and $1 \leq j < l \leq p.$
Asymptotically in $d$ these are $k {{p(p-1)} \over {2}}$ equations, bilinear in the entries $s^j.$ This leads to ${{p(p-1)} \over {2}}$ polynomials
of degree $2d,$ thus to $\sim k(p-1)pd$ conditions (the vanishing of the corresponding coefficients). These are the conditions which force
a section $s: C \to C \times G$ to map to $\Sigma.$ To make sure that $s$ maps to $Z$ we need further conditions.
The number of linear conditions describing the relative linear space $F$ is asympotically equal to
$$ dim \Sigma - \dim F = p(m-p) - k{{p(p-1)} \over {2}} - n+1.$$
Since the Pl\"ucker coordinates are polynomials of degree $p$ in the $s^j,$ these linear conditions lead to the vanishing of
$p(m-p) - k{{p(p-1)} \over {2}} -n+1 $ polynomials in $x$ (coordinate on $\bP_1)$), each of degree $pd$, and so we get
$$\sim (p(m-p) -k {{p(p-1)} \over {2}} -n+1)pd$$
conditions. \\
In total we obtain asymptotically
$$ kp(p-1)d + \big( p(m-p) - k{{p(p-1)} \over {2}} - n+1\big) pd $$
conditions, while we have $mpd$ parameters to our disposal. Hence we get a section as soon as we have more parameters than conditions,
and this leads immediately to inequality (1) in (4.3). \\
Notice that we should also have checked that $s^1 \wedge \ldots \wedge s^p \ne 0$. But this is done in exactly the same way as in (4.2). \\

{\it Case 2:} $l = p-1, k = 1.$ Here the arguments are completely analogous. We again have $mpd$ parameters. The number of conditions forcing
$s$ to be in $\Sigma$ is asymptotically $dpl = p(p-1)d$, because the dimension of $(p-1)-$linear alternating forms on $\comp ^p$ is
$p = {{p} \choose {p-1}}$, and each coefficient of the relative form is a polynomial of degree $p-1$ in the $s^j,$ which leads
to a polynomial of degree $p(p-1)d$ in $x.$ \\
The number of conditions forcing $s$ to be in $F$ is as in case (1)
$$ (\dim \sigma - n)pd = \big( p(m-p) - p - n +1\big) pd, $$
and by the inequality (2) we have at least as many parameters as equations so that $s$ exists. \qed

{\bf (4.5)} We are discussing the case $g = 7.$ Here $m = 9, p = 4$ whereas $l = 2; k = 1.$ In that case the inequalities are not satisfied,
so there is no conclusion.
However we can still conclude that $\Sigma \to B$ has at least two sections (actually a 1-dimensional family), say $s_1 $ and $s_2.$ We then choose a base
$e_i$ of $\comp ^9$ carefully so that the equations of $\Sigma \subset C \times G \subset C \times \bP^N$ take a nice form. \\
Using these two sections, we choose the $e_i$ such that $s_1$ is the subspace generated by $e_1, \ldots e_4 $ and $q(e_i,e_j) = \delta_{ij}.$
Moreover we may choose $e_6, \ldots, e_9$ such that $s_2$ is given by the $4-$plane defined by those $e_i.$ We can also achieve that
$q(e_i,e_{5+j}) = 0$ and $q(e_{5+i},e_{5+j}) = \delta_{ij}.$ Finally let $e_5$ be in the orthogonal complement of the $e_i$ already chosen.
Denoting by $x$ again the variable on $C$ and denoting by $s_{i,j} $ the Pl\"ucker coordinates of  the section $s$ ($1 \leq i \leq 4$; $1 \leq j
\leq 9)$,
we can write the equations for a section $s$ to be in $Z$ in the following form:
$$ a_i(x)s^2_{i,5} + b_i(t)s_{i,5+i} = 0 $$
for $1 \leq i \leq 4 $ and
$$ a_{ij}(x)s^2_{i,5} + c_{ij}(x)s_{i,5+j} + d_{ji}(x)s_{j,5+i} $$
for the obvious set of $i,j$ (these are ${{4} \choose {2}} = 6$ equations). If we associate to $s_{i,5}$ the weight ${d} \over {2} $ and to all
other coordinates weight $d$, then we obtain, adding the $7$ additional linear equation for describing $Z \subset \Sigma$, in total
$10d + 7d = 17d$ conditions. On the other hand we have $16d + 4{d \over 2} = 18 d$ parameters, which is bigger than the number of
conditions. This settles the case $g = 7.$  \\
\\ It remains to treat Fano fibrations whose fibers have $b_2 \geq 2.$ This is rather tedious and works by running through
the classification of Mori and Mukai. Of course the geometry of the Mori contractions plays the decisive role.
Here is just one argument as example - details will be worked out in another paper.
Suppose that some (hence every) smooth fiber $X_c$ admits only one birational contraction. Then we can contract fiberwise
to obtain $X \rightharpoonup X' \to C.$ If the birational contraction of $X_c$ contracts the divisor to a point, then
$X' $ has a section. Otherwise we get a distinguished curve $C_c \subset X'_c,$ i.e. we obtain a ruled surface
$S \subset X' \to C,$ which has a section.

\vspace{1cm}
\small
\begin{tabular}{ll}
Fr\'ed\'eric Campana &Thomas Peternell \\
D\'epartement de math\'ematiques &
Mathematisches Institut  \\
Universit\'e de Nancy I & Universit\"at Bayreuth \\
BP 239, 54506 Vandoeuvre les Nancy Cedex & 95440 Bayreuth \\
France & Germany\\
frederic.campana@antares.iecn.u-nancy.fr & thomas.peternell@uni-bayreuth.de   \\
\end{tabular}

\vspace{1cm}

\begin{tabular}{ll}
Aleksandr Pukhlikov & \\
Steklov Mathematical Institute & Mathematisches Institut \\
Gubkina 8 & Universit\"at Bayreuth \\
117966 Moscow & 95440 Bayreuth \\
Russia & Germany \\
pukh@mi.ras.ru &pukh@btm8x5.mat.uni-bayreuth.de \\
\end{tabular}


\newpage




\begin{thebibliography}{999999}

\bibitem[Do80]{Do80} Dolgachev,I.: Newton polyhedra and factorial rings. J. pure and appl. algebra 18, 253-258 (1980)

\bibitem[Do81]{Do81} Dolgachev,I.: Weighted projective spaces. Lecture Notes in Math. 956, 34-71. Springer 1981

\bibitem[GHS01]{GHS01} Graber,T.;Harris,J.;Starr,J.: Families of rationally connected varieties. Preprint, june 2001

\bibitem[IP99]{IP99} Iskovskikh,V.A.;Prokhorov,Yu.G.: Fano varieties. Encyclopaedia in Math.Sci.47. Springer 1999

\bibitem[Ts36]{Ts36} Tsen,C.: Quasi-algebraisch abgeschlossene Funktionenk\"orper. J.Chinese Math.1, 81-92 (1936)

\bibitem[Ko96]{Ko96} Koll\' ar J.: Rational curves on algebraic varieties. Springer 1996

\bibitem[KMM92]{KMM92} Koll\' ar J., Miyaoka,Y.  and Mori S.: Rationally connected varieties. J.Alg.Geom.1, 429-448 (1992)

\end{thebibliography}
\end{document}